\documentclass[10pt,notitlepage]{amsart}
\usepackage{amsmath} 
\usepackage{amscd}
\usepackage{amssymb}
\usepackage{enumerate}
\usepackage{indentfirst}

\font\smallrm=cmr8

\renewcommand{\:}{\colon}

\newcommand{\wt}{\widetilde}
\newcommand{\w}{\omega}
\newcommand{\wh}{\widehat}
\newcommand{\ox}{\otimes}
\newcommand{\cC}{{\mathcal C}}
\newcommand{\cU}{{\mathcal U}}

\newcommand{\cF}{{\mathcal F}}

\newcommand{\cO}{{\mathcal O}}
\newcommand{\cT}{{\mathcal T}}
\newcommand{\cA}{{\mathcal A}}
\newcommand{\cS}{{\mathcal S}}
\newcommand{\cK}{{\mathcal K}}
\newcommand{\cR}{{\mathcal R}}

\newcommand{\cJ}{{\mathcal J}}
\newcommand{\cX}{{\mathcal X}}

\newcommand{\cL}{{\mathcal L}}
\newcommand{\CC}{{\text{\bf C}}}
\newcommand{\QQ}{{\text{\bf Q}}}

\newtheorem{theorem}{Theorem}[section]

\newtheorem{corollary}[theorem]{Corollary}
\newtheorem{proposition}[theorem]{Proposition}

\theoremstyle{definition} 
\newtheorem{definition}[theorem]{Definition}
\newtheorem{example}[theorem]{Example}

\newtheorem{subsec}[theorem]{Construction}
\theoremstyle{plain}

\begin{document}

\

\title[\smallrm Jets of singular foliations]
{Jets of singular foliations}
\author[\smallrm Eduardo Esteves]{Eduardo Esteves}

\thanks{
2000 Mathematics Subject Classification: 
14F10, 13N15, 37F75.}
\thanks{Supported by CNPq, Proc.~478625/03-0 
and 301117/04-7, 
FAPERJ, Proc. E-26/170.418/00, and 
CNPq/FAPERJ, Proc. E-26/171.174/03.}

\begin{abstract} Given a singular foliation 
satisfying locally everywhere the Frobenius 
condition, even at the singularities, 
we show how to construct its global sheaves of jets. 
Our construction is purely formal, 
and thus applicable in a variety of contexts. 
\end{abstract}

\maketitle

\section{Introduction}

Let $M$ be a complex manifold of complex dimension 
$m$. A holomorphic foliation $\cL$ of 
dimension $n$ of $M$ is a decomposition of $M$ in 
complex submanifolds, called leaves, of 
dimension $n$. Also, locally the leaves must pile up nicely, 
like the fibers of a holomorphic map. In other words, 
for each point $p$ of $M$ there must exist an open 
neighborhood $U$ and a holomorphic submersion 
$\varphi\: U\to V$ to an open subset 
$V\subseteq \CC^{m-n}$ such that the fibers of 
$\varphi$ are the intersection of the leaves with 
$U$. We say that $\varphi$ defines $\cL$ on 
$U$. 

A holomorphic foliation $\cL$ of $M$ induces a vector 
subbundle of the tangent bundle $T_M$ of $M$: for 
each $p\in M$, the vector subspace of $T_{M,p}$ is the 
tangent space at $p$ of the unique leaf passing 
through $p$. Thinking in dual terms, $\cL$ 
induces a surjection $w\:T^*_M\to E$ from the bundle 
of holomorphic 1-forms to a holomorphic rank-$n$ 
vector bundle $E$. The bundle $E$, also denoted by 
$T^*_\cL$, is regarded as the bundle of 1-forms 
of $\cL$.

Not all surjections $w\:T^*_M\to E$ to a 
holomorphic vector bundle $E$ arise from foliations. 
The necessary and sufficient condition for this is 
given by the Frobenius Theorem: locally at each point 
$p$ of $M$, choose a trivialization of $E$, and 
consider the vector fields $X_1,\dots,X_n$ induced by $w$; 
if their Lie brackets $[X_i,X_j]$ can be expressed as 
sums $\sum_\ell g_\ell X_\ell$, where the $g_\ell$ 
are holomorphic functions on a neighborhood of $p$, 
then $w$ arises from a foliation.

The surjection $w$ can be seen, locally on 
an open subset $U\subseteq M$ for which there is a 
submersion $\varphi\:U\to V\subseteq\CC^{m-n}$ 
defining $\cL$, as the natural map 
$T^*_U\to T^*_{U/V}$ from the bundle of 1-forms on 
$U$ to the bundle of relative 1-forms on $U$ over 
$V$. Also, on such a $U$, we may consider the natural 
surjection $J^q_U\to J^q_{U/V}$ from the bundle 
$J^q_U$ of $q$-jets (or principal parts of order $q$) 
on $U$ to the bundle $J^q_{U/V}$ of relative jets of 
$\varphi$, for 
each integer $q\geq 0$. 
These patch to form surjections 
$w^q\:J^q_M\to J^q_\cL$ to bundles $J^q_\cL$ that can 
be regarded as the bundles of $q$-jets of the 
foliation. 

But what happens if all the data are algebraic?
More precisely, assume $M$ is algebraic, $E$ and $w$ are 
algebraic, and there is an algebraic trivialization 
of $E$ at each $p\in M$ such that 
the resulting vector fields $X_1,\dots,X_n$ 
are involutive, i.e. satisfy 
$[X_i,X_j]=\sum_\ell g^{i,j}_\ell X_\ell$ 
for $g_\ell^{i,j}$ algebraic. Are the bundles 
$J^m_\cL$ algebraic? In 
principle, they are just holomorphic, since the 
local submersions $\varphi$ from which they arise 
are just holomorphic, constructed by means of the 
implicit function theorem.

Moreover, it is rare for a projective manifold to admit 
interesting foliations. For this reason, 
one has started to study \emph{singular foliations}, in a 
variety of ways. For instance, 
a singular foliation of $M$ of dimension $n$ may be defined 
to be a map $w\:T^*_M\to E$ to a rank-$n$ holomorphic 
vector bundle $E$ which, on a dense open subset 
$M^0\subseteq M$, arises from a foliation $\cL$. We 
still regard $E$ as the bundle of 1-forms of the 
foliation. But the bundles of jets $J^m_\cL$ are, in 
principle, only defined on $M^0$. 
Under which conditions do they extend to $M$? 

In the present paper, we will show that if all the 
data are algebraic, then the bundles $J^m_\cL$ are 
algebraic. Furthermore, if the same Frobenius' 
conditions, appropriately formulated, are verified 
at each point of $M-M^0$, then the bundles $J^m_\cL$ 
extend to the whole $M$. For the 
proofs, we will completely bypass Frobenius Theorem, 
giving an entirely formal construction of the bundles 
of jets that applies in many categories, for instance 
that of differentiable manifolds, or of schemes over 
any base. Furthermore, not only will we consider maps 
to bundles $E$, but also to sheaves of modules, 
locally free or not, obtaining thus sheaves of jets.

Our construction of the sheaves of jets is by 
``iteration'', so will only apply in characteristic 
zero. For an approach in positive 
characteristic, albeit limited in scope, 
see \cite{E1} or \cite{LT3}. 
From now on, all rings are assumed to be 
$\QQ$-algebras.

Here is what we do. 
Let $X$ be a topological space, $\cO_B$ a sheaf of 
$\QQ$-algebras 
and $\cO_X$ a sheaf of $\cO_B$-algebras. Let $\cF$ be 
a sheaf of $\cO_X$-modules and $D\:\cO_X\to\cF$ an 
$\cO_B$-derivation. We may think of $\cO_X$ as the 
sheaf of ``functions'' on $X$ and of $\cO_B$ as the 
sheaf of ``constant functions.'' For each 
integer $i\geq 0$, 
let $\cT^i(\cF)$ be the tensor product over $\cO_X$ 
of $i$ copies of $\cF$, and denote by $\cS^i(\cF)$ 
and $\cA^i(\cF)$ its symmetric and exterior 
quotients. Let $\cS(\cF)$ be the direct product of 
all the $\cS_i(\cF)$ for all integers 
$i\geq 0$, with its natural graded 
$\cO_X$-algebra structure. We may think of 
$\cS(\cF)$ as the sheaf of ``formal power series'' 
on the sections of $\cF$.

We define a $D$-connection to be a map of $\cO_B$-modules 
$\gamma\:\cF\to\cT^2(\cF)$ satisfying 
         $$\gamma(am)=D(a)\ox m+a\gamma(m)$$
for all local sections $a$ of $\cO_X$ and $m$ of 
$\cF$. We will see in Construction \ref{extcon} 
how $\gamma$ can be used to iterate $D$ to 
obtain a sequence of maps 
$T_i:\cF\to\cS^i(\cF)\ox\cF$, for $i=0,1,\dots$, 
where $T_0:=\text{id}_{\cF}$, $T_1=\gamma$, and 
the $T_i$ satisfy properties similar to that of 
a connection, i.e. Equations~\ref{connection}. We 
call any sequence of maps $T=(T_0,T_1,\dots)$ with 
these properties, for any $D$-connection $\gamma$, 
an extended $D$-connection.

From an extended $D$-connection $T$ we get a map 
$h\:\cO_X\to\cS(\cF)$ of 
$\cO_B$-algebras, with $h_0:=\text{id}_{\cO_X}$ and 
$h_1:=D$, by letting $h_i(a)$ be the class in 
$\cS^i(\cF)$ of $(1/i)T_{i-1}D(a)$ for each local 
section $a$ of $\cO_X$; see Construction~\ref{existhas}.
We may think of $h$ as a way of computing 
``Taylor series'' of functions on $X$. In this way, 
$\cS(\cF)$ may be regarded as a sheaf of jets. We call 
such an $h$ an iterated Hasse derivation. 

However, $D$-connections are usually local gadgets. So, 
to be able to patch the local maps $h$, we need 
$D$-connections 
to be canonical. But there is nothing canonical 
about $\gamma$: for every $\cO_X$-linear map 
$\nu\:\cF\to\cT^2(\cF)$, the sum $\gamma+\nu$ is also 
a $D$-connection! So we study special $D$-connections.

A $D$-connection $\gamma$ is 
called flat if the image of $\gamma D(\cO_X)$ in 
$\cA^2(\cF)$ is zero. When such $\gamma$ exists, 
we say that $D$ is integrable. Our main result, 
Theorem~\ref{existint}, shows that, given a flat 
$D$-connection $\gamma$, there is an extended 
$D$-connection $T$ with $T_1=\gamma$, which is also flat, 
meaning that $T_iD(\cO_X)$ lies in the subsheaf of 
$\cS^i(\cF)\ox\cF$ generated locally by
       $$\sum_{\ell=1}^im_1m_2\cdots m_{\ell-1}m_{\ell+1}
         m_{\ell+2}\cdots m_i\ox m_\ell$$
for all local sections $m_1,\dots,m_i$ of $\cF$, for each 
$i\geq 1$.

Now, our Proposition \ref{iterated} says that a flat, extended 
$D$-connection is ``comparable'' to any other extended 
$D$-connection. So, Theorem \ref{existint} and 
Proposition \ref{iterated} can be coupled to yield that all 
iterated Hasse derivations are equivalent; see 
Corollary \ref{equivint}. More precisely, if $h$ and $h'$ are 
iterated Hasse derivations, there is an 
$\cO_X$-algebra automorphism 
$\phi\:\cS(\cF)\to\cS(\cF)$ such that $h'=\phi h$. 
Furthermore, the degree-$i$ part of this $\phi$ is 
zero if $i<0$ and the identity if $i=0$.

Now, for the patching we also need the 
automorphisms $\phi$ to be 
canonical, so that cocycle conditions are satisfied. 
For this, we make a technical assumption on $D$, 
that bounds its singularities, and holds in all 
applications we know of; see Corollary~\ref{equivint} and 
the remark thereafter. 

Finally, assume that $D$ is locally 
integrable, and has bounded singularities, in the sense 
alluded to above. The patching of the local iterated 
Hasse derivations and the $\cO_X$-algebra 
automorphisms comparing 
them is straightforward. We obtain an 
$\cO_X$-algebra $\cJ$ and a map $h\:\cO_X\to\cJ$ of 
$\cO_B$-algebras. Also, since the local 
$\cO_X$-algebra automorphisms do not decrease 
degrees, and their degree-0 parts are the identity, 
$\cJ$ comes naturally with a filtration by 
$\cO_X$-algebra quotients $\cJ^q$, for $q=0,1,\dots$, 
and natural exact sequences
       	$$0\to\cS^q(\cF)\to\cJ^q\to\cJ^{q-1}\to0$$
for each $q>0$; see Construction~\ref{jets}.
We say that $\cJ^q$ is the sheaf of $q$-jets of $D$.

How does this formal construction fit with the 
geometric setting? If $M$ is a holomorphic manifold, 
let $(X,\cO_X)$ be the ringed space where $X$ is the 
underlying topological space of $M$ and 
$\cO_X$ is its sheaf of 
holomorphic functions. Let $\cO_B$ be 
the sheaf of constant complex 
functions. A map of vector bundles $w\:T^*_M\to E$ 
corresponds to a derivation $D\:\cO_X\to\cF$, where 
$\cF$ is the sheaf of holomorphic sections of $E$. 
To say that $w$ arises from a foliation on an open 
subset $U\subseteq M$ is equivalent to say that 
$\cF|_U=\cO_UD(\cO_U)$ and $D|_U$ is locally 
integrable on $U$; see Example~\ref{loccon}. 
Now, assume that  
$D$ is locally integrable on the whole $X$, and 
that $D(\cO_X)$ generates $\cF$ as an $\cO_X$-module 
on a dense open subset $M^0\subseteq M$. 
Then there exists a sheaf of $q$-jets $\cJ^q$ on 
$X$, as explained above.  Also, $w$ defines a 
foliation $\cL$ on $M^0$, and $\cJ^q|_{M^0}$ is 
the sheaf of sections 
of the bundle of $q$-jets $J^q_\cL$;
see Example~\ref{jetsfol}.

Bundles of jets associated to foliations or 
derivations were considered by a number of people. 
In algebraic geometry, to my knowledge, the first 
was Letterio Gatto, who in his thesis 
\cite{G1} constructed jets from a family 
of stable curves $f\:\cX\to\cS$ and its canonical 
derivation $\cO_X\to\omega_{X/S}$, where $\omega_{X/S}$ 
is the relative dualizing 
sheaf. Afterwards, in \cite{E1}, jets were 
constructed for more general 
families, of local complete intersection curves, over 
any base and in any characteristic.

Also, Dan Laksov and Anders Thorup constructed  
bundles of jets in a series of 
articles in different setups; see \cite{LT1},\cite{LT2} 
and \cite{LT3}. In characteristic zero, 
their more general work is \cite{LT2}. Actually, 
in \cite{LT2}, Laksov and Thorup construct larger 
sheaves of ``jets'', that are naturally 
noncommutative. The true 
generalization of the sheaf of jets is what they call 
``symmetric'' jets. They show that the sheaf of 
(symmetric) jets is uniquely defined when $\cF$ is free 
and has an $\cO_X$-basis under which $D$ can be 
expressed using commuting derivations of $\cO_X$. As 
we have observed above, the uniqueness of the definition is 
important for patching. However, the commutativity is 
stronger than Frobenius' conditions, at least in the 
algebraic category --- in the analytic category, 
at nonsingular points, the local 
existence of commuting derivations follows from the 
existence of the foliation. The present work arose from the 
feeling that the Frobenius' conditions should be enough to 
construct sheaves of jets. 

There have already 
been applications of the sheaves of jets 
associated to a foliation or a derivation. 
They were used by Gatto 
\cite{G2}, and Gatto and myself \cite{EG} in 
enumerative aplications. They were used by myself in 
understanding limits of ramification points 
\cite{E2}, and of Weierstrass points, with Nivaldo 
Medeiros, \cite{EM1} and \cite{EM2}. 
They were also used by Jorge Vit\'orio Pereira in 
the study of foliations of the projective space \cite{P}. 

Finally, we will see an example where the 
integrability condition holds and $\cF$ is not 
locally free; see Example~\ref{nonG}. 
That will be the example of the canonical derivation 
on a special non-Gorenstein 
curve, arguably the 
simplest non-Gorenstein unibranch curve there is, 
whose complete local ring at the singular point is 
of the form $\CC[[t^3,t^4,t^5]]$, as a subring 
of $\CC[[t]]$. It could be that the 
integrability condition holds for the canonical derivation 
on any curve, Gorenstein or not. If so, the sheaf of 
jets might be used to define Weierstrass points on 
non-Gorenstein integral curves, and show that these 
points are limits of Weierstrass points on nearby 
curves, in the way done by Robert Lax and 
Carl Widland for Gorenstein curves; 
see \cite{LW} and \cite{W}. 
However, this problem will not be pursued 
here.

Here is a layout of the article. In Section 2, we define 
connections, extended connections and iterated Hasse 
derivations, and explain a few preliminary constructions. 
In Section 3, we define integrable derivations and 
flat (extended) connections, and show that a flat, extended 
connection is comparable with 
any other extended connection. 
Finally, in Section 4, we show that, if a derivation $D$ 
is integrable, then flat, extended connections 
exist, and all iterated Hasse derivations 
are equivalent; then we construct the sheaves 
of jets for locally integrable derivations.

This work started as a joint work with Letterio Gatto. 
However, he felt he did not contribute to it 
as much as he wished. 
Though a few discussions with him were vital to how this 
work came to be, and though I feel that this could be 
classified as a joint work, I had to respect his decision 
to not coauthor it. 
Anyway, being the only thing left for me to do, I thank him 
for his great contributions.

\section{Derivations and connections}

\begin{subsec}\setcounter{equation}{0} 
(\emph{Tensor operations}) Let $X$ be a topological space and 
$\cO_X$ a sheaf of $\QQ$-algebras. 
Let  $\cF$ be a sheaf of $\cO_X$-modules. Denote by 
	$$\cT(\cF):=\prod_{n=0}^\infty\cT^n(\cF),\quad
	\cS(\cF):=\prod_{n=0}^\infty\cS^n(\cF),\quad
	\cA(\cF):=\prod_{n=0}^\infty\cA^n(\cF)$$
the \emph{formal} \emph{tensor}, \emph{symmetric} and 
\emph{exterior} graded sheaf of $\cO_X$-algebras of $\cF$, 
respectively. (Note that we take the direct product and not 
the direct sum.)

Set $\cR^0(\cF):=\cO_X$. Also, set 
$\cR^n(\cF):=\cS^{n-1}(\cF)\ox\cF$ 
for each integer $n\geq 1$, and
        $$\cR(\cF):=\prod_{n=0}^{\infty}\cR^n(\cF).$$
Then $\cR(\cF)$ is a graded 
$\cO_X$-algebra quotient of $\cT(\cF)$, in a natural way.

As usual, we let $\cT_+(\cF)$, $\cS_+(\cF)$, $\cA_+(\cF)$ and 
$\cR_+(\cF)$ denote the ideals generated by 
formal sums with zero constant terms in each of the 
indicated $\cO_X$-algebras. 

We view $\cT(\cF)$ as a sheaf of algebras, with the 
(noncommutative) product induced by tensor product. 
So, given local sections 
$m_1,\dots,m_n$ of $\cF$, we let $m_1\cdots m_n$ denote their 
product in $\cT^n(\cF)$. Also, we view $\cS(\cF)$, $\cA(\cF)$ 
and $\cR(\cF)$ as sheaves of $\cT(\cF)$-algebras, and 
$\cS(\cF)$ as a sheaf of $\cR(\cF)$-algebras, under the 
natural quotient maps. So, given a local section $\omega$ of 
$\cT^n(\cF)$ (resp. $\cR^n(\cF)$), we will 
use the same symbol $\omega$ to denote its image in 
$\cS^n(\cF)$, $\cA^n(\cF)$ or $\cR^n(\cF)$ 
(resp. $\cS^n(\cF)$). These simplifications should not lead to 
confusion, and will clean the notation enormously. 

Define the \emph{switch operator} 
$\sigma\:\cR_+(\cF)\to\cR_+(\cF)$ as 
the homogeneous $\cO_X$-linear map of degree 0 given by 
$\sigma|_{\cF}:=0$, and on each $\cR^n(\cF)$, for $n\geq 2$, 
by the formula:
	$$\sigma(m_1\cdots m_n)=\sum_{i=1}^{n-1}m_nm_{n-1}
        \cdots 
	m_{n-i+2}m_{n-i+1}m_1m_2\cdots m_{n-i-1}m_{n-i}$$
for all local sections $m_1,\dots,m_n$ of $\cF$. The reader may check that $\sigma$ is actually well-defined on 
$\cR_+(\cF)$, and not only on $\cT_+(\cF)$. 

Let $\sigma^\star:=1+\sigma$. Notice that 
$\sigma^\star$ factors through $\cS_+(\cF)$. For each integer 
$n\geq 1$, let $\cK^n(\cF):=\sigma^\star(\cR^n(\cF))$. 
Then $\cK^n(\cF)$ is also the kernel of $n-\sigma^\star$. In 
particular, $\cK^2(\cF)$ is the kernel of the surjection 
$\cT^2(\cF)\to\cA^2(\cF)$. Indeed, that $\cK^n(\cF)$ is in the 
kernel of $n-\sigma^\star$ follows from the equality 
        $$\sigma^\star\sigma^\star|_{\cR^n(\cF)}=
        n\sigma^\star|_{\cR^n(\cF)},$$
a fact checked locally. And if $\w$ is a local section of 
$\cR^n(\cF)$ such that $(n-\sigma^\star)(\w)=0$, then 
$\w=\sigma^\star((1/n)\w)$, and thus $\w$ is a local 
section of $\cK^n(\cF)$. 

Put
        $$\cK_+(\cF):=\prod_{n=1}^\infty\cK^n(\cF).$$
\end{subsec}

\begin{definition} \setcounter{equation}{0} 
Let $X$ be a topological space, $\cO_B$ a 
sheaf of $\QQ$-algebras and 
$\cO_X$ a sheaf of $\cO_B$-algebras. 
Let  $\cF$ be a sheaf of $\cO_X$-modules. 
A \emph{Hasse derivation of $\cF$} is a 
map of $\cO_B$-algebras 
	$$h=(h_0,h_1,\dots)\:\cO_X\longrightarrow
        \prod_{i=0}^\infty\cS^i(\cF)$$
with $h_0=\text{id}_{\cO_X}$. Let $D\:\cO_X\to\cF$ be an 
$\cO_B$-derivation. We say that $h$ \emph{extends $D$} 
if $h_1=D$.
\end{definition}

If $h=(h_0,h_1,\dots)$ is a Hasse derivation of $\cF$, then 
$h_1\:\cO_X\to\cF$ is an $\cO_B$-derivation of 
$\cF$. Conversely, given an $\cO_B$-derivation 
$D\:\cO_X\to\cF$, we may ask when there is a Hasse derivation 
$h=(h_0,h_1,\dots)$ of $\cF$ extending $D$. 
We will see in Construction~\ref{existhas} 
that such $h$ exists when there is a $D$-connection.
 
\begin{definition} \setcounter{equation}{0} 
Let $X$ be a topological space, $\cO_B$ a 
sheaf of $\QQ$-algebras and 
$\cO_X$ a sheaf of $\cO_B$-algebras. 
Let $\cF$ be a sheaf of $\cO_X$-modules and 
$D\:\cO_X\to\cF$ an $\cO_B$-derivation. 
A \emph{$D$-connection} is a map of $\cO_B$-modules
	$$\gamma\:\cF\to\cT^2(\cF)$$
satisfying
	\begin{equation}\label{eq:03}
	\gamma(am)=D(a)m+a\gamma(m)
	\end{equation}
for each local sections $a$ of $\cO_X$ and $m$ of $\cF$.
\end{definition}

\begin{example}\label{nodal}\setcounter{equation}{0} 
There may not exist a $D$-connection. For 
instance, let $X$ be the union of two transversal lines in the 
plane, or $X:=\text{Spec}(\CC[x,y]/(xy))$. Let 
$\Omega^1_X$ be the sheaf of differentials, and 
$D\:\cO_X\to\Omega^1_X$ the universal $\CC$-derivation. The sheaf 
$\Omega^1_X$ is generated by $D(x)$ and $D(y)$, and the sheaf 
of relations is generated by the single relation $yD(x)+xD(y)=0$. 
In particular, $D(x)$ and $D(y)$ are $\CC$-linearly independent at 
the node. Suppose there were a $D$-connection 
$\gamma\:\Omega^1_X\to\cT^2(\Omega^1_X)$. Then
     $$0=\gamma(yD(x)+xD(y))
     =D(y)D(x)+D(x)D(y)+y\gamma(D(x))+
     x\gamma(D(y)).$$
However, $D(x)D(x)$, $D(x)D(y)$, $D(y)D(x)$ and 
$D(y)D(y)$ are linearly independent sections of 
$\cT^2(\Omega^1_X)$ at the 
node. Hence the above relation is not possible. 
\end{example}

When a $D$-connection $\gamma$ exists, it is not unique, since 
for every $\cO_X$-linear map $\lambda\:\cF\to\cT^2(\cF)$, 
the sum $\gamma+\lambda$ is a $D$-connection. However, 
these are all the $D$-connections.

A $D$-connection allows us to iterate $D$ 
to a Hasse derivation, 
as we will explain in Construction \ref{existhas}. 
First, we will see how to extend a connection.

\begin{subsec}\label{extcon}\setcounter{equation}{0} 
(\emph{Extending connections}) 
Let $X$ be a topological space, $\cO_B$ a 
sheaf of $\QQ$-algebras and 
$\cO_X$ a sheaf of $\cO_B$-algebras. 
Let $\cF$ be a sheaf of $\cO_X$-modules and 
$D\:\cO_X\to\cF$ an 
$\cO_B$-derivation. Let $\gamma\:\cF\to\cT^2(\cF)$ be a 
$D$-connection. 
Define a homogeneous map of degree 1 of $\cO_B$-modules,
	\begin{equation}\label{eq:01}
	\nabla\:\cR_+(\cF)\to\cR_+(\cF),
	\end{equation}
given on each graded part $\cR^n(\cF)$ by
	\begin{equation}\label{eq:02}
	\nabla(m_1\cdots m_n):=\sum_{i=1}^n 
	m_1\cdots m_{i-1}\gamma(m_i)m_{i+1}\cdots m_n
	\end{equation}
for each local sections $m_1,\dots,m_n$ of $\cF$.

At first, it seems $\nabla$ would 
be a well-defined map from $\cT_+(\cF)$ 
to $\cT_+(\cF)$. This would indeed 
be true, were $\gamma$ a map of 
$\cO_X$-modules. But $\gamma$ is not! To check that $\nabla$, 
as given above, is well-defined, 
we need to check the following 
three properties for all local sections 
$m_1,\dots,m_i,m'_i,\dots,m'_n$ of $\cF$ and 
$a$ of $\cO_X$, and each 
permutation $\tau$ of $\{1,\dots,n-1\}$:
	\begin{enumerate}
	\item For each $i=1,\dots,n$,
	  $$\nabla(m_1\cdots(m_i+m'_i)\cdots m_n)=
	  \nabla(m_1\cdots m_i\cdots m_n)
	  +\nabla(m_1\cdots m'_i\cdots m_n).$$
	\item For each $i,j=1,\dots,n$,
	  $$\nabla(m_1\cdots(am_i)\cdots m_n)=
	  \nabla(m_1\cdots(am_j)\cdots m_n).$$ 
	\item $\nabla(m_1\cdots m_{n-1}m_n)=
	\nabla(m_{\tau(1)}\cdots m_{\tau(n-1)}m_n)$.
	\end{enumerate}
The first (multilinearity) and third (symmetry) 
properties are left for the reader to check. 
The second property is the key to why $\nabla$ 
must take values in $\cR_+(\cF)$, so let us check it: 
from the definition of $\nabla$, and using 
$\gamma(am_i)=D(a)m_i+a\gamma(m_i)$, we get
	$$\nabla(m_1\cdots(am_i)\cdots m_n)=
        m_1\cdots m_{i-1}D(a)m_i\cdots m_n+
	a\sum_{j=1}^nm_1\cdots\gamma(m_j)\cdots m_n.$$
So $\nabla(m_1\cdots(am_i)\cdots m_n)$ 
would depend on $i$, were 
$\nabla$ to take values in $\cT_+(\cF)$. But instead, $\nabla$ 
takes values in $\cR_+(\cF)$, and hence
	$$m_1\cdots m_{i-1}D(a)m_i\cdots m_n=
        D(a)m_1\cdots m_{i-1}m_i\cdots m_n.$$
So the second property (scalar multiplication) 
is checked, and the 
three properties imply that $\nabla$ is well-defined. Also, we 
proved the formula
	$$\nabla(a\w)=D(a)\w+a\nabla(\w)$$
for all local sections $\w$ of 
$\cR_+(\cF)$ and $a$ of $\cO_X$. 

Now, for each integer $n\geq 1$, put
        $$T_n:=\frac{1}{n!}\nabla^n|_{\cF}\:
        \cF\to\cR^{n+1}(\cF).$$
Also, set $T_0:=\text{id}_\cF$. 
Then $T_n=(1/n)\nabla T_{n-1}$ for 
each $n\geq 1$. Furthermore, for each integer $i\geq 1$, and each 
local sections $a$ of $\cO_X$ and $m$ of $\cF$,
        \begin{equation}\label{connection}
	T_i(am)=aT_i(m)+\sum_{j=1}^i\frac{1}{j}
	T_{j-1}D(a)T_{i-j}(m).
	\end{equation}
Indeed, Formula \eqref{connection} holds for $i=1$, because 
$T_1=\gamma$ and $\gamma$ is a $D$-connection. And if, 
by induction, Formula \eqref{connection} holds for  
a certain $i\geq 1$, then 
        \begin{align*}
	T_{i+1}(am)=&\frac{1}{(i+1)}\nabla T_i(am)\\
	=&\frac{1}{(i+1)}\nabla\Big(aT_i(m)+
	\sum_{j=1}^i\frac{1}{j}T_{j-1}D(a)T_{i-j}(m)\Big)\\
	=&\frac{1}{(i+1)}\Big(a\nabla T_i(m)+D(a)T_i(m)\Big)\\
	+&\frac{1}{(i+1)}\sum_{j=1}^i\frac{1}{j}\Big(\nabla
	T_{j-1}D(a)T_{i-j}(m)+T_{j-1}D(a)\nabla T_{i-j}(m)\Big)\\
	=&aT_{i+1}(m)+\frac{1}{(i+1)}D(a)T_i(m)\\
	+&\frac{1}{(i+1)}\sum_{j=1}^i\Big(T_jD(a)T_{i-j}(m)+
	\frac{(i+1-j)}{j}T_{j-1}D(a)T_{i+1-j}(m)\Big)\\
	=&aT_{i+1}(m)+\frac{1}{(i+1)}D(a)T_i(m)\\
	+&\sum_{j=2}^i\Big(
	\frac{1}{j}T_{j-1}D(a)T_{i+1-j}(m)\Big)+
	\frac{1}{(i+1)}T_iD(a)m+\frac{i}{(i+1)}D(a)T_i(m)\\
	=&aT_{i+1}(m)+\sum_{j=1}^{i+1}
	\frac{1}{j}T_{j-1}D(a)T_{i+1-j}(m).
	\end{align*}
The maps $T_n$ form an extended $D$-connection, according to the 
definition below.
\end{subsec}
	 
\begin{definition}\setcounter{equation}{0} 
Let $X$ be a topological space, $\cO_B$ a 
sheaf of $\QQ$-algebras and 
$\cO_X$ a sheaf of $\cO_B$-algebras. 
Let $\cF$ be a sheaf of $\cO_X$-modules 
and $D\:\cO_X\to\cF$ an 
$\cO_B$-derivation. Let $n$ 
be a positive integer or $n:=\infty$. A map of 
$\cO_B$-modules
	$$T=(T_0,T_1,T_2,\dots)\:\cF\longrightarrow
        \prod_{i=0}^n\cR^{i+1}(\cF)$$ 
is called an {\em extended $D$-connection} if $T_0=\text{id}_\cF$ 
and Formula \eqref{connection} holds for each $i\geq 1$ and 
all local sections $a$ of $\cO_X$ and $m$ of $\cF$.
\end{definition}

If nothing is noted otherwise, extended 
$D$-connections are assumed {\em full}, 
that is, with $n:=\infty$.

\begin{subsec}\label{existhas}\setcounter{equation}{0}
(\emph{Hasse derivations arising from $D$-connections}) 
Let $X$ be a topological space, $\cO_B$ a 
sheaf of $\QQ$-algebras and 
$\cO_X$ a sheaf of $\cO_B$-algebras. 
Let $\cF$ be a sheaf of $\cO_X$-modules, $D\:\cO_X\to\cF$ an 
$\cO_B$-derivation and $T\:\cF\to\cR_+(\cF)$ an extended 
$D$-connection. Notice that the map 
$T_1\:\cF\to\cF\ox\cF$ is a 
$D$-connection. However, $T$ need not arise from $T_1$ as in 
Construction~\ref{extcon}.

Put $h_0:=\text{id}_{\cO_X}$, and for 
each integer $i\geq 1$ let $h_i\:\cO_X\to\cS^i(\cF)$ be the 
$\cO_B$-linear map given by $h_i(a):=(1/i)T_{i-1}D(a)$ for each 
local section $a$ of $\cO_X$. Then
        $$h:=(h_0,h_1,h_2,\dots)\:\cO_X\longrightarrow
        \prod_{i=0}^{\infty}\cS^i(\cF)$$
is a Hasse derivation of $\cF$ extending $D$. Indeed, 
clearly $h_1=D$. Now, if $a$ and $b$ are local sections of 
$\cO_X$, and $i\geq 1$, then
        \begin{align*}
	h_i(ab)=&\frac{1}{i}T_{i-1}\Big(aD(b)+bD(a)\Big)\\
	=&\frac{1}{i}\Big(aT_{i-1}D(b)+\sum_{j=1}^{i-1}\frac{1}{j}
	T_{j-1}D(a)T_{i-1-j}D(b)\Big)\\
	+&\frac{1}{i}\Big(bT_{i-1}D(a)+\sum_{j=1}^{i-1}
	\frac{1}{j}T_{j-1}D(b)T_{i-1-j}D(a)\Big)\\
	=&ah_i(b)+bh_i(a)\\
	+&\frac{1}{i}\sum_{j=1}^{i-1}\Big(
	\frac{1}{j}T_{j-1}D(a)T_{i-1-j}D(b)+\frac{1}{i-j}
	T_{i-1-j}D(b)T_{j-1}D(a)\Big)\\
	=&ah_i(b)+bh_i(a)+\sum_{j=1}^{i-1}\Big(\frac{1}{j(i-j)}
	T_{j-1}D(a)T_{i-1-j}D(b)\Big)\\
	=&\sum_{j=0}^ih_j(a)h_{i-j}(b),
	\end{align*}
where in the fourth equality we used that the computation is done 
in $\cS^i(\cF)$.
\end{subsec}

\begin{definition} Let $X$ be a topological space, $\cO_B$ a 
sheaf of $\QQ$-algebras and 
$\cO_X$ a sheaf of $\cO_B$-algebras. 
Let $\cF$ be a sheaf of $\cO_X$-modules, $D\:\cO_X\to\cF$ an 
$\cO_B$-derivation and $h$ a Hasse derivation extending $D$. 
We say that $h$ is \emph{iterated} if there is an 
extended $D$-connection $T$ such that $h_i=(1/i)T_{i-1}D$ 
for each $i\geq 1$.
\end{definition}

\section{Flat connections and integrable derivations}

\begin{definition}\setcounter{equation}{0} 
Let $X$ be a topological space, $\cO_B$ a 
sheaf of $\QQ$-algebras and $\cO_X$ 
a sheaf of $\cO_B$-algebras. 
Let $\cF$ be a sheaf of $\cO_X$-modules and 
$D\:\cO_X\to\cF$ an $\cO_B$-derivation. 
We say that a $D$-connection $\gamma\:\cF\to\cT^2(\cF)$ is 
\emph{flat} if $\gamma D(\cO_X)\subseteq\cK^2(\cF)$. 
We say that $D$ is \emph{integrable} if there exists 
a flat $D$-connection. 
\end{definition}

\begin{example}\label{loccon}\setcounter{equation}{0} 
Let $X$ be a topological space, $\cO_B$ a sheaf of 
$\QQ$-algebras and $\cO_X$ a sheaf of $\cO_B$-algebras. 
Let $\cF$ be 
the free sheaf of $\cO_X$-modules with basis $e_1,\dots,e_n$. 
Let $D\:\cO_X\to\cF$ be an $\cO_B$-derivation. Then 
$D=D_1e_1+\cdots+D_ne_n$, where the $D_i$ are 
$\cO_B$-derivations of $\cO_X$. Conversely, a 
$n$-tuple $(D_1,\dots,D_n)$ of $\cO_B$-derivations of 
$\cO_X$ defines an $\cO_B$-derivation of $\cF$. 

There is a natural $D$-connection 
$\gamma\:\cF\to\cT^2(\cF)$, satisfying
	$$\gamma(\sum_{i=1}^na_ie_i)=
        \sum_{i=1}^n\sum_{j=1}^nD_j(a_i)e_je_i$$
for all local sections $a_1,\dots,a_n$ of $\cO_X$. 
Any other $D$-connection is of the form $\gamma+\nu$, 
where $\nu\:\cF\to\cT^2(\cF)$ is a map of 
$\cO_X$-modules. The map $\nu$ is defined by global 
sections $c^{j,i}_\ell$ of $\cO_X$, for 
$1\leq i,j,\ell\leq n$, satisfying
	$$\nu(e_\ell)=
        \sum_{i=1}^n\sum_{j=1}^nc^{j,i}_\ell e_je_i.$$
To say that $(\gamma+\nu)D(a)=0$ in $\cA^2(\cF)$ for 
a local section $a$ of $\cO_X$ is to say that
	$$\sum_{i=1}^n\sum_{j=1}^nD_jD_i(a)e_je_i+
        \sum_{\ell=1}^n\sum_{i=1}^n\sum_{j=1}^n
        c^{j,i}_\ell D_\ell(a)e_je_i=0$$
in $\cA^2(\cF)$ or, equivalently, 
	$$D_jD_i(a)-D_iD_j(a)=\sum_{\ell=1}^n
	(c^{i,j}_\ell-c^{j,i}_\ell)D_\ell(a)$$
for all distinct $i$ and $j$. In other words, $D$ is 
integrable if and only if the collection 
$\{D_1,\dots,D_n\}$ is \emph{involutive}, i.e., if and 
only if there are sections $b^{j,i}_\ell$ of $\cO_X$ 
such that
	$$[D_j,D_i]=\sum_{\ell=1}^n
        b^{j,i}_\ell D_\ell$$
for all distinct $i$ and $j$.

Let $T$ be the extended $D$-connection of 
Construction \ref{extcon}, derived from $\gamma$. 
From the definition of $\gamma$ we have 
$T_q(e_i)=0$ for each integer $q>0$ and 
each $i=1,\dots,n$. Suppose $\gamma D=0$ in $\cA^2(\cF)$, 
or in other words $[D_j,D_i]=0$ for all $i$ and $j$. Then 
the Hasse derivation $h$ associated to $T$ satisfies
        \begin{equation}\label{hcom}
	h_q(a)=\sum_{j_1=1}^n\cdots\sum_{j_q=1}^n
	\frac{D_{j_1}\cdots D_{j_q}(a)}{q!}
	e_{j_1}\cdots e_{j_q}
	\end{equation}
for each integer $q>0$ and each local section $a$ 
of $\cF$.
\end{example}

\begin{example}\label{nonG}\setcounter{equation}{0} 
It is not necessary that $\cF$ be locally 
free for a flat $D$-connection to exist. For instance, let 
$X:=\text{Spec}(\CC[t^3,t^4,t^5])$. Viewed as a sheaf 
of regular meromorphic differentials, Rosenlicht-style, 
the dualizing sheaf $\omega_X$ is generated by $dt/t^2$ 
and $dt/t^3$. 
Let $\eta_1:=dt/t^2$ and $\eta_2:=dt/t^3$. The following 
relations generate all relations $\eta_1$ and 
$\eta_2$ satisfy with coefficients in $\cO_X$:
       \begin{equation}\label{releta}
	 t^3\eta_1=t^4\eta_2,\quad
	 t^4\eta_1=t^5\eta_2,\quad 
	 t^5\eta_1=t^6\eta_2.
       \end{equation}
The composition 
of the universal derivation with the canonical map 
$\Omega^1_X\to\omega_X$ yields a derivation 
$D\:\cO_X\to\omega_X$ satisfying
     $$D(t^3)=3t^4\eta_1=3t^5\eta_2,\quad
       D(t^4)=4t^5\eta_1=4t^6\eta_2,\quad
       D(t^5)=5t^6\eta_1=5t^7\eta_2.$$
So $D(a)\eta_i=0$ in $\cA^2(\omega_X)$ for each local 
section $a$ of $\cO_X$ and each $i=1,2$. 

Define $\gamma\:\omega_X\to\cT^2(\omega_X)$ by letting
     $$\gamma(a\eta_1+b\eta_2)=D(a)\eta_1+D(b)\eta_2+
     4t^3a\eta_2\eta_2+3b\eta_1\eta_1$$
for all local sections $a$ and $b$ of $\cO_X$.
To check that $\gamma$ is well defined we need only check 
that the values of $\gamma$ on both sides of the three 
relations \eqref{releta} agree. This is 
the case; for instance, 
      $$\gamma(t^3\eta_1)=3t^4\eta_1\eta_1+4t^6\eta_2\eta_2
      =7t^4\eta_1\eta_1=4t^5\eta_1\eta_2+3t^4\eta_1\eta_1=
      \gamma(t^4\eta_2).$$
Now, since $D(a)\eta_1=D(b)\eta_2=0$ in $\cA^2(\omega_X)$, 
we have that $\gamma=0$ in $\cA^2(\omega_X)$. 
So $\gamma$ is a flat $D$-connection, and hence $D$ is 
integrable.
\end{example}

\begin{definition}\setcounter{equation}{0} 
Let $X$ be a topological space, $\cO_B$ a 
sheaf of $\QQ$-algebras and 
$\cO_X$ a sheaf of $\cO_B$-algebras. 
Let $\cF$ be 
a sheaf of $\cO_X$-modules, and let $D\:\cO_X\to\cF$ be an 
$\cO_B$-derivation. An extended $D$-connection 
$T\:\cF\to\cR_+(\cF)$ is called \emph{flat} if 
$TD(\cO_X)\subseteq\cK_+(\cF)$.
\end{definition}

We will see in Theorem \ref{existint} that flat extended 
$D$-connections exist, when $D$ is integrable. Also, by 
Proposition \ref{iterated} below, any two of them are 
``comparable''. First, a piece of notation.

\begin{subsec}\setcounter{equation}{0} (\emph{Generating maps}) 
Let $X$ be a topological space, 
$\cO_X$ a sheaf of $\QQ$-algebras and $\cF$ a sheaf of 
$\cO_X$-modules. Let $n$ be a 
positive integer or $n:=\infty$. Let 
        $$\lambda=(\lambda_0,\lambda_1,\dots)\:
        \cF\longrightarrow\prod_{i=0}^n\cR^{i+1}(\cF)$$
be a map of $\cO_X$-modules. For each integer $p>0$ and 
each sequence $i_1,\dots,i_p$ of nonnegative indices at 
most equal to $n$, let $q:=i_1+\cdots+i_p$ and define 
	$$(\lambda_{i_1}\cdots\lambda_{i_p})\:
         \cT^p(\cF)\to\cR^{p+q}(\cF)$$
to be the map of $\cO_X$-modules satisfying 
	$$(\lambda_{i_1}\cdots\lambda_{i_p})
        (m_1\cdots m_p):=
        \lambda_{i_1}(m_1)\cdots\lambda_{i_p}(m_p)$$
for all local sections $m_1,\dots,m_p$ of $\cF$.

The maps $(\lambda_{i_1}\cdots\lambda_{i_p})$ are not 
defined on $\cR^p(\cF)$, but the sum 
	$$s_q(\lambda):=\sum_{i_1+\cdots+i_p=q}
	(i_p+1)\big(\lambda_{i_1}\cdots
	\lambda_{i_p}\big)\:\cR^p(\cF)
	\longrightarrow\cR^{p+q}(\cF)$$ 
is, for all integers $p>0$ and each integer $q$ with 
$0\leq q\leq n$. Analogously, the sum
        $$\wt s_q(\lambda):=\sum_{i_1+\cdots+i_p=q}
	\big(\lambda_{i_1}\cdots\lambda_{i_p}\big)\:\cS^p(\cF)
	\longrightarrow\cS^{p+q}(\cF)$$
is well-defined, for all integers $p>0$ 
and each integer $q$ with 
$0\leq q\leq n$. 

Notice that, for each local section $\w$ of $\cK^p(\cF)$,
        \begin{equation}\label{swts}
	  s_q(\lambda)(\w)=\frac{p+q}{p}\wt s_q(\lambda)(\w) 
	  \text{ in $\cS^{p+q}(\cF)$}.
	  \end{equation}
Indeed, locally,
        $$\w=\sum_{i=1}^pm_1\cdots\wh{m_i}\cdots m_pm_i$$
for local sections $m_1,\dots,m_p$ of $\cF$. Thus, in $\cS^{p+q}(\cF)$,
	\begin{align*}
	s_q(\lambda)(\w)=&
	\sum_{i=1}^p\sum_{j_1+\cdots+j_p=q}
	(j_p+1)\prod_{s=1}^{i-1}\lambda_{j_s}(m_s)
	\prod_{s=i}^{p-1}
	\lambda_{j_s}(m_{s+1})\lambda_{j_p}(m_i)\\
	=&\sum_{j_1+\cdots+j_p=q}\sum_{i=1}^p
	(j_i+1)\prod_{s=1}^p\lambda_{j_s}(m_s)\\
	=&(p+q)\sum_{j_1+\cdots+j_p=q}
	\prod_{s=1}^p\lambda_{j_s}(m_s)\\
	=&\frac{p+q}{p}\sum_{j_1+\cdots+j_p=q}
	\sum_{i=1}^p\prod_{s=1}^{i-1}\lambda_{j_s}(m_s)
	\prod_{s=i}^{p-1}\lambda_{j_s}(m_{s+1})
	\lambda_{j_p}(m_i)\\
	=&\frac{p+q}{p}\sum_{j_1+\cdots+j_p=q}
	(\lambda_{j_1}\cdots\lambda_{j_p})(\w)\\
	=&\frac{p+q}{p}\wt s_q(\lambda)(\w).
	\end{align*}
\end{subsec}

\begin{proposition}\label{iterated}\setcounter{equation}{0} 
Let $X$ be a topological space, $\cO_B$ a 
sheaf of $\QQ$-algebras and 
$\cO_X$ a sheaf of $\cO_B$-algebras. 
Let $\cF$ be a sheaf of $\cO_X$-modules 
and $D\:\cO_X\to\cF$ an 
$\cO_B$-derivation. Let $n$ be a positive integer. Let 
        \begin{align*}
	T=(T_0,T_1,\dots,T_n)\:&\cF\longrightarrow
	\prod_{i=0}^n\cR^{i+1}(\cF),\\
	S=(S_0,S_1,\dots,S_n)\:&\cF\longrightarrow
	\prod_{i=0}^n\cR^{i+1}(\cF)
	\end{align*}
be two $\cO_B$-linear maps. Assume that 
$T_iD(\cO_X)\subseteq\cK^{i+1}(\cF)$ 
for each $i=0,\dots,n-1$. 
Then any two of the following three 
statements imply the third:
\begin{enumerate}
\item The map $S$ is an extended $D$-connection.
\item The map $T$ is an extended $D$-connection.
\item There is a (unique) map of $\cO_X$-modules 
  $$\lambda=(\lambda_0,\dots,\lambda_n)\:\cF
  \longrightarrow\prod_{i=0}^n\cR^{i+1}(\cF)$$ 
such that $\lambda_0=\text{\rm id}_{\cF}$ and
	\begin{equation}\label{equiviter}
	S_i=\sum_{\ell=0}^is_{i-\ell}(\lambda)T_\ell
	\end{equation}
for each $i=0,1,\dots,n$.
\end{enumerate}
\end{proposition}

\begin{proof} 
We will argue by induction. For $n=1$, the 
proposition simply says that, given a $D$-connection $T_1$, 
a map of $\cO_B$-modules $S_1\:\cF\to\cR^2(\cF)$ is a 
$D$-connection if and only if $S_1-T_1$ is $\cO_X$-linear, 
a fact already observed.

Suppose now that $n\geq 2$, and that the statement of the 
proposition is known for $n-1$ in place of $n$. So we may 
assume that $(S_0,\dots,S_{n-1})$ and 
$(T_0,\dots,T_{n-1})$ are extended 
$D$-connections, and that there exists a map of 
$\cO_X$-modules 
       $$\lambda=(\lambda_0,\dots,\lambda_{n-1})\:\cF
       \longrightarrow\prod_{i=0}^{n-1}\cR^{i+1}(\cF)$$ 
such that $\lambda_0=\text{id}_{\cF}$ and 
Equations \eqref{equiviter} hold 
for each $i<n$. 

Define 
        $$R_n:=\sum_{\ell=1}^{n-1}
        s_{n-\ell}(\lambda)T_\ell.$$
We need only show that, for each local sections 
$a$ of $\cO_X$ and $m$ of $\cF$,
        \begin{equation}\label{RST}
        R_n(am)-aR_n(m)+
        \sum_{z=0}^{n-1}\frac{T_{n-z-1}D(a)}{n-z}T_z(m)=
        \sum_{z=0}^{n-1}\frac{S_{n-z-1}D(a)}{n-z}S_z(m).
	\end{equation}
Indeed, if $S$ and $T$ are $D$-connections, Formula \eqref{RST} 
implies that $S_n-R_n-T_n$ is 
$\cO_X$-linear. So, setting 
$\lambda_n:=(1/(n+1))(S_n-R_n-T_n)$, 
Equation \eqref{equiviter} holds for $i=n$ as well. 
Conversely, if there is an $\cO_X$-linear 
map $\lambda_n\:\cF\to\cR^{n+1}(\cF)$ such that Equation 
\eqref{equiviter} holds for $i=n$, then 
$(n+1)\lambda_n+R_n=S_n-T_n$. So, from Formula \eqref{RST} 
we see that $S$ is a $D$-connection if and only if $T$ is. 

Now, on the one hand, since $(T_0,\dots,T_{n-1})$ is an 
extended $D$-connection, 
	\begin{align*}
	R_n(am)-aR_n(m)=&\sum_{\ell=1}^{n-1}
	\sum_{j_0+\cdots+j_\ell=n-\ell}
	(j_\ell+1)\big(\lambda_{j_0}\cdots
	\lambda_{j_\ell}\big)
	\Big(T_\ell(am)-aT_\ell(m)\Big)\\
	=&\sum_{\ell=1}^{n-1}
	\sum_{j_0+\cdots+j_\ell=n-\ell}
	\sum_{p=0}^{\ell-1}\frac{j_\ell+1}{\ell-p}
	\big(\lambda_{j_0}\cdots\lambda_{j_\ell}\big)
	\Big(T_{\ell-1-p}D(a)T_p(m)\Big).
	\end{align*}
Thus the left-hand side of Formula \eqref{RST} is equal to
        \begin{equation}\label{lhsRST}
	\sum_{\ell=1}^n\sum_{j_0+\cdots+j_\ell=n-\ell}
	\sum_{p=0}^{\ell-1}\frac{j_\ell+1}{\ell-p}
	\big(\lambda_{j_0}\cdots\lambda_{j_\ell}\big)
	\Big(T_{\ell-1-p}D(a)T_p(m)\Big).
	\end{equation}

On the other hand, using Equations \eqref{equiviter} 
for $i<n$, the right-hand side of \eqref{RST} becomes
	$$\sum_{z=0}^{n-1}\frac{1}{n-z}
	\bigg(\sum_{k=0}^{n-z-1}
	\sum_{j_0+\cdots+j_k=n-z-1-k}(j_k+1)\big(
	\lambda_{j_0}\cdots\lambda_{j_k}\big)T_kD(a)
	\bigg)\bigg(\sum_{p=0}^z\w_{z-p}\bigg),$$
where
        $$\w_\ell:=\sum_{j'_0+\cdots+j'_p=\ell}
	(j'_p+1)\big(\lambda_{j'_0}\cdots
	\lambda_{j'_p}\big)T_p(m)$$
for each $\ell=0,\dots,n-1$. 
Now, for each $z=0,\dots,n-1$ and $k=0,\dots,n-z-1$, 
using that $T_kD(a)$ is a local section of 
$\cK^{k+1}(\cF)$, Formula \eqref{swts} yields the following equation in 
$\cS^{n-z}(\cF)$:
	$$\sum(j_k+1)
	(\lambda_{j_0}\cdots\lambda_{j_k})T_kD(a)
	=\sum\frac{n-z}{k+1}
	(\lambda_{j_0}\cdots\lambda_{j_k})T_kD(a),$$
where the sum on both sides runs over the $(k+1)$-tuples 
$(j_0,\dots,j_k)$ such that $j_0+\cdots+j_k=n-z-1-k$. 
Thus, introducing 
$\ell:=k+p+1$, the right-hand side of \eqref{RST} becomes
	$$\sum_{\ell=1}^n\sum_{p=0}^{\ell-1}\frac{1}{\ell-p}
	\sum_{z=p}^{n-\ell+p}
	\sum_{j_0+\cdots+j_{\ell-p-1}=n-z-\ell+p}
	(\lambda_{j_0}\cdots\lambda_{j_{\ell-p-1}})
	T_{\ell-p-1}D(a)\w_{z-p},$$
whence, introducing $u:=z-p$, equal to
	$$\sum_{\ell=1}^n\sum_{p=0}^{\ell-1}\frac{1}{\ell-p}
	\sum_{u=0}^{n-\ell}
	\sum_{j_0+\cdots+j_{\ell-p-1}=n-\ell-u}
	(\lambda_{j_0}\cdots\lambda_{j_{\ell-p-1}})
	T_{\ell-p-1}D(a)\w_u,$$
which is equal to \eqref{lhsRST}.
\end{proof}

\section{Jets}

\begin{theorem}\label{existint}
\setcounter{equation}{0}
Let $X$ be a topological space, $\cO_B$ a 
sheaf of $\QQ$-algebras and 
$\cO_X$ a sheaf of $\cO_B$-algebras. 
Let $\cF$ be a sheaf of $\cO_X$-modules and 
$D\:\cO_X\to\cF$ an 
$\cO_B$-derivation. 
If $D$ is integrable, then there exists a 
flat, extended $D$-connection. 
\end{theorem}

\begin{proof} Since $D$ is 
integrable, there is a flat $D$-connection 
$T_1\:\cF\to\cT^2(\cF)$. Set $T_0:=\text{id}_{\cF}$. 
Suppose, by induction, that for an integer $n\geq 2$ 
we have constructed an extended $D$-connection 
	$$T=(T_0,T_1,\dots,T_{n-1})\:\cF\longrightarrow
	\prod_{i=0}^{n-1}\cR^{i+1}(\cF).$$
We will also suppose the maps 
$T_i$ satisfy one additional property, Equations~\eqref{Ti}, 
after we make a definition.

For each $j=0,\dots,n-1$, define a map of $\cO_B$-modules 
$T'_j\:\cR^2(\cF)\to\cR^{j+2}(\cF)$ by letting 
	$$T'_j(m_1m_2):=\sum_{i=0}^jT_i(m_1)T_{j-i}(m_2)$$
for all local sections $m_1$ and $m_2$ of $\cF$. 
To check that $T'_j$ is well defined, we need only check 
that 
$\sum_iT_i(am_1)T_{j-i}(m_2)=\sum_iT_i(m_1)T_{j-i}(am_2)$ 
for each local section $a$ of $\cO_X$. 
In fact, using \eqref{connection}, we see that both sides 
are equal to 
	$$\sum_{i=0}^jaT_i(m_1)T_{j-i}(m_2)+
        \sum_{\ell=1}^{j}\sum_{i=0}^{j-\ell}
	\frac{1}{\ell}T_{\ell-1}D(a)
	T_i(m_1)T_{j-\ell-i}(m_2).$$
Furthermore, we see from this computation that
	$$T'_j(a\w)=aT'_j(\w)+\sum_{i=1}^j
	\frac{1}{j+1-i}T_{j-i}D(a)T'_{i-1}(\w)$$
for all local sections $\w$ of $\cR^2(\cF)$ 
and $a$ of $\cO_X$.

Also, notice that $\sigma^\star T'_{i-1}(1-\sigma)=0$. 
Indeed, for local sections $m_1$ and $m_2$ of $\cF$,
	\begin{align*}
	\sigma^\star T'_{i-1}(1-\sigma)(m_1m_2)=&
	\sigma^\star T'_{i-1}(m_1m_2-m_2m_1)\\
	=&\sigma^\star\bigg(\sum_{j=0}^{i-1}
	\Big(T_j(m_1)T_{i-1-j}(m_2)-T_{i-1-j}(m_2)T_j(m_1)\Big)\bigg)\\
	=&\sum_{j=0}^{i-1}\Big(\sigma^\star
	\big(T_j(m_1)T_{i-1-j}(m_2)\big)-
	\sigma^\star\big(T_{i-1-j}(m_2)T_j(m_1)\big)\Big),
	\end{align*}
which is zero because 
$\sigma^\star(\w_1\w_2)=\sigma^\star(\w_2\w_1)$ for 
all local sections $\w_1$ and $\w_2$ of $\cR_+(\cF)$. 
Then 
	\begin{equation}\label{Tprime}
	T'_{i-1}(1-\sigma)=
        \frac{i-\sigma}{i+1}T'_{i-1}(1-\sigma).
	\end{equation}

Now, suppose that
	\begin{equation}\label{Ti}
	(i-\sigma)\big((i+1)T_i-T'_{i-1}(1-\sigma)T_1\big)=0
	\end{equation}
for each $i=1,\dots,n-1$. (Notice that 
Equation \eqref{Ti} holds automatically for 
$i=1$, because $(1-\sigma)\sigma^\star(\cR^2(\cF))=0$.) 
Also, from Equation \eqref{Ti}, and the flatness of $T_1$, 
we get $T_iD(\cO_X)\subseteq\cK^{i+1}(\cF)$ 
for each $i=1,\dots,n-1$. 

Let
	$$T:=\frac{1}{n+1}T'_{n-1}(1-\sigma)T_1.$$
Then
	\begin{equation}\label{Tcon}
	(n-\sigma)\big(T(am)-aT(m)-\sum_{i=1}^n
	\frac{1}{i}T_{i-1}D(a)T_{n-i}(m)\big)=0
	\end{equation}
for all local sections $m$ of $\cF$ and $a$ of $\cO_X$. 
Indeed, 
	\begin{align*}
	(n-\sigma)T(am)=&\frac{n-\sigma}{n+1}
	T'_{n-1}(1-\sigma)T_1(am)\\
	=&\frac{n-\sigma}{n+1}
	\bigg(T'_{n-1}\Big(D(a)m-mD(a)+a(1-\sigma)T_1(m)\Big)\bigg)\\
	=&\frac{n-\sigma}{n+1}\bigg(T_{n-1}D(a)m+
	\sum_{i=1}^{n-1}T_{n-1-i}D(a)T_i(m)\bigg)\\
	-&\frac{n-\sigma}{n+1}\bigg(mT_{n-1}D(a)+
	\sum_{i=1}^{n-1}T_i(m)T_{n-1-i}D(a)\bigg)\\
	+&\frac{n-\sigma}{n+1}\bigg(aT'_{n-1}(1-\sigma)T_1(m)\bigg)\\
	+&\frac{n-\sigma}{n+1}\bigg(\sum_{i=1}^{n-1}
	T_{n-i-1}D(a)T'_{i-1}\frac{(1-\sigma)}{n-i}T_1(m)\bigg).
	\end{align*}
Now, first observe that, for each $i=0,\dots,n-1$, we have 
$(n-i-\sigma^\star)T_{n-1-i}D(a)=0$, since 
$T_{n-1-i}D(a)$ is a local section of $\cK^{n-i}(\cF)$. Then
        \begin{align*}
	(n-\sigma)\Big(T_i(m)T_{n-1-i}D(a)\Big)=&
        (n-\sigma)\Big(T_i(m)\frac{\sigma^\star}{n-i}T_{n-1-i}D(a)\Big)\\
	=&-\frac{n-\sigma}{n-i}\Big(T_{n-1-i}D(a)\sigma^\star T_i(m)\Big).
	\end{align*}
Also, from \eqref{Tprime} and \eqref{Ti},
	\begin{align*}
	T_{n-i-1}D(a)T'_{i-1}(1-\sigma)T_1(m)=&
	T_{n-i-1}D(a)\frac{i-\sigma}{i+1}T'_{i-1}
	(1-\sigma)T_1(m)\\
	=&T_{n-i-1}D(a)(i-\sigma)T_i(m)
	\end{align*}
for each $i=1,\dots,n-1$. So
        $$\frac{n-\sigma}{n+1}\bigg(
	T_{n-1}D(a)m-mT_{n-1}D(a)\bigg)=\frac{n-\sigma}{n}
	\bigg(T_{n-1}D(a)m\bigg),$$
and, for each $i=1,\dots,n-1$, 
	$$\frac{n-\sigma}{n+1}\bigg(
	T_{n-1-i}D(a)T_i(m)-T_i(m)T_{n-1-i}D(a)+
	T_{n-i-1}D(a)T'_{i-1}\frac{1-\sigma}{n-i}T_1(m)\bigg)$$
is equal to 
        $$\frac{n-\sigma}{n+1}\bigg(T_{n-1-i}D(a)
	\big(1+\frac{1+\sigma}{n-i}+\frac{i-\sigma}{n-i}\big)
	T_i(m)\bigg),$$
whence equal to 
        $$\frac{n-\sigma}{n-i}\bigg(T_{n-1-i}D(a)T_i(m)\bigg).$$
Applying these equalities in the above expression for 
$(n-\sigma)T(am)$ we get Equation~\eqref{Tcon}.

We want to show that there exists a map of $\cO_B$-modules 
$T_n\:\cF\to\cR^{n+1}(\cF)$ such that \eqref{Ti} 
holds for $i=n$, and 
such that $(T_0,\dots,T_n)$ is an extended $D$-connection. 
First we claim that there exists a map of $\cO_B$-modules 
$T_n\:\cF\to\cR^{n+1}(\cF)$ such that 
$(T_0,\dots,T_n)$ is an extended $D$-connection. 
Indeed, from $T_1$ construct an extended $D$-connection 
$(S_0,\dots,S_n)$ by iteration, as described in 
Construction \ref{extcon}. 
By Proposition \ref{iterated}, there is a map 
of $\cO_X$-modules 
	$$\lambda=(\lambda_0,\lambda_1,\dots,
        \lambda_{n-1})\:\cF\longrightarrow
	\prod_{i=0}^{n-1}\cR^{i+1}(\cF)$$ 
such that $\lambda_0=\text{id}_{\cF}$ and such that 
Equations \eqref{equiviter} hold for $i=0,\dots,n-1$. Now, 
just set $\lambda_n:=0$ in Equation \eqref{equiviter} 
for $i=n$, and let it define $T_n$. Then, by 
Proposition~\ref{iterated}, the 
map $(T_0,\dots,T_n)$ is an extended $D$-connection.

The above map $T_n$ does not necessarily make 
\eqref{Ti} hold for $i=n$. So, rename it by $U$. At 
any rate, since $(T_0,\dots,T_{n-1},U)$ is an extended 
$D$-connection, it follows from Equation \eqref{Tcon} 
that $(n-\sigma)(U-T)$ is $\cO_X$-linear. Set
	$$T_n:=U-\frac{(n-\sigma)}{n+1}(U-T).$$
Then $T_n$ differs from $U$ by an $\cO_X$-linear map, 
and thus $(T_0,\dots,T_n)$ is an extended $D$-connection. 
Now,
	$$(n-\sigma)T_n=(n-\sigma)U-
        \frac{(n-\sigma)^2}{n+1}(U-T).$$
Since 
	$$(n-\sigma)^2|_{\cR^{n+1}(\cF)}=
        (n+1)(n-\sigma)|_{\cR^{n+1}(\cF)},$$
we get $(n-\sigma)T_n=(n-\sigma)T$. So \eqref{Ti} holds for $i=n$.

The induction argument is complete, showing that there is an 
infinite extended $D$-connection 
	$$T=(T_0,T_1,\dots)\:\cF\longrightarrow
        \prod_{i=0}^\infty \cR^{i+1}(\cF)$$
such that \eqref{Ti} holds for each $i\geq 1$, and thus 
$T_iD(\cO_X)\subseteq\cK^{i+1}(\cF)$ 
for each $i\geq 0$.
\end{proof}

\begin{definition}\setcounter{equation}{0} 
Let $X$ be a topological space, $\cO_B$ a 
sheaf of $\QQ$-algebras and $\cO_X$ a sheaf 
of $\cO_B$-algebras. Let  
$\cF$ be a sheaf of $\cO_X$-modules. 
Two Hasse derivations $h$ and $h'$ 
of $\cF$ are said to be 
\emph{equivalent} if there is an $\cO_X$-algebra 
automorphism $\phi$ of $\cS(\cF)$ such that 
$\phi_0|_{\cF}=\text{id}_{\cF}$ and $h'=\phi h$. 
We say that $h$ and $h'$ are 
\emph{canonically equivalent} when there is only 
one such automorphism.
\end{definition} 

\begin{corollary}\label{equivint}
\setcounter{equation}{0} 
Let $X$ be a topological space, $\cO_B$ a 
sheaf of $\QQ$-algebras and $\cO_X$ a 
sheaf of $\cO_B$-algebras. 
Let  $\cF$ be a sheaf of $\cO_X$-modules and $D\:\cO_X\to\cF$  
an $\cO_B$-derivation. 
Let $h$ and $h'$ be iterated Hasse 
derivations of $\cF$ extending $D$. If $D$ is 
integrable, then $h$ and $h'$ are equivalent. 
Furthermore, if $\nu D(\cO_X)\neq 0$ for 
every nonzero $\cO_X$-linear map 
$\nu\:\cF\to\cS(\cF)$, 
then $h$ and $h'$ are canonically equivalent. 
\end{corollary}

\begin{proof} First, we prove the existence of an 
equivalence. By Theorem~\ref{existint}, there is a 
flat, extended $D$-connection 
$T=(T_0,T_1,\dots)$. 
We may suppose 
$h$ arises from $T$. Let $S=(S_0,S_1,\dots)$ be an 
extended $D$-connection from which $h'$ arises.

By Proposition~\ref{iterated}, 
there is a map of $\cO_X$-modules 
	$$\lambda=(\lambda_0,\lambda_1,\dots)\:\cF
        \longrightarrow
	\prod_{i=0}^\infty\cR^{i+1}(\cF)$$ 
such that $\lambda_0=\text{id}_{\cF}$ and
	\begin{equation}\label{rel}
	S_i=\sum_{\ell=0}^is_{i-\ell}(\lambda)T_\ell
	\quad\text{for each $i\geq 0$.}
\end{equation}

Let $\phi\:\cS(\cF)\to\cS(\cF)$ be the 
map of $\cO_X$-algebras whose 
graded part $\phi_q$ of degree $q$ satisfies 
$\phi_q|_{\cS^p(\cF)}=\wt s_q(\lambda)|_{\cS^p(\cF)}$ for all 
integers $p>0$ if $q\geq 0$, and $\phi_q=0$ if $q<0$.
Since $\lambda_0=\text{id}_{\cF}$, the homogeneous 
degree-0 part $\phi_0$ is the identity, and thus 
$\phi$ is an automorphism. We claim that 
$h'=\phi h$.

Indeed, clearly $h'_0=(\phi h)_0$. Now, since $T$ is 
flat, 
$T_\ell D(\cO_X)\subseteq\cK^{\ell+1}(\cF)$ 
for each $\ell\geq 0$. Thus, for each $i\geq 0$ and each local section 
$a$ of $\cO_X$, using Equations \eqref{swts} and \eqref{rel}, the 
following equalities hold on $\cS^{i+1}(\cF)$:
	$$h'_{i+1}(a)=\frac{S_iD(a)}{i+1}=
        \sum_{\ell=0}^is_{i-\ell}(\lambda)\frac{T_\ell D(a)}{i+1}
	=\sum_{\ell=0}^i\wt s_{i-\ell}(\lambda)\frac{T_\ell D(a)}{\ell+1}
	=\sum_{\ell=0}^i\phi_{i-\ell}h_{\ell+1}(a).$$
Since $\phi_{i+1}|_{\cO_X}=0$, we have $h'_{i+1}=(\phi h)_{i+1}$. 
So $h'=\phi h$.

Now, assume that $\nu D(\cO_X)\neq 0$ for 
every nonzero $\cO_X$-linear map 
$\nu\:\cF\to\cS(\cF)$. 
Let $\phi$ be an $\cO_X$-algebra automorphism of 
$\cS(\cF)$ such that $\phi_0|_{\cF}=\text{id}_{\cF}$ 
and $h'=\phi h$. To see that $\phi$ is unique, we 
just need to show that $\phi_q|_{\cF}$ is uniquely 
defined for each $q\geq 0$. We do it by induction. 
Since $h'=\phi h$, for each $q\geq 0$ the following 
equality holds:
         $$h'_{q+1}=\phi_qh_1+\phi_{q-1}h_2+\cdots+
         \phi_1h_q+h_{q+1}.$$
Then $\phi_qD$ is determined by 
$\phi_1,\dots,\phi_{q-1}$. Since $\phi_q$ is 
$\cO_X$-linear and 
takes values in $S^{q+1}(\cF)$, it follows from 
our extra assumption that $\phi_q$ is determined.
\end{proof} 

If $(X,\cO_X)$ is a Noetherian scheme over a 
Noetherian $\QQ$-scheme $(B,\cO_B)$, and $\cF$ is 
a locally free sheaf on $X$, then, for 
every $\cO_B$-derivation 
$D\:\cO_X\to\cF$ such that $\cF$ is generated by 
$D(\cO_X)$ at the associated points of $X$, we have 
$\nu D(\cO_X)\neq 0$ for all nonzero 
$\cO_X$-linear maps $\nu\:\cF\to\cS(\cF)$.

\begin{subsec}\label{jets}
\setcounter{equation}{0} (\emph{Jets}) 
Let $X$ be a topological space, $\cO_B$ a 
sheaf of $\QQ$-algebras, 
and $\cO_X$ a sheaf of $\cO_B$-algebras. 
Let $\cF$ be a sheaf of $\cO_X$-modules and $D\:\cO_X\to\cF$  
an $\cO_B$-derivation. Assume that $D\:\cO_X\to\cF$ 
is locally integrable, and that the sheaf of 
$\cO_X$-linear maps from $\cF$ to $\cS(\cF)$ 
sending $D(\cO_X)$ to zero has only trivial 
local sections. Let $\cU$ be the 
collection of open subspaces $U\subseteq X$ such 
that $D|_U$ is integrable. For each $U\in\cU$, 
there exist iterated Hasse derivations 
extending $D|_U$. Let $\cC_U$ be the 
collection of these Hasse derivations. 
By Corollary \ref{equivint}, for any two 
$h,h'\in\cC_U$ there is a unique $\cO_U$-algebra 
automorphism $\phi_{h,h'}$ of $\cS(\cF)|_U$ 
such that $h'=\phi_{h,h'}h$. 
Now, consider the collection of all 
the $h\in\cC_U$ for all $U\in\cU$. 
Consider also the collection of all the 
$\phi_{h,h'}$ for all $U\in\cU$ and all 
$h,h'\in\cC_U$. If $h,h',h''\in\cC_U$, then 
$h''=\phi_{h',h''}\phi_{h,h'}h$. 
From the uniqueness of $\phi_{h,h''}$, we get 
$\phi_{h',h''}\phi_{h,h'}=\phi_{h,h''}$. The cocycle 
condition being satisfied, the 
$\phi_{h,h'}$ patch the $\cS(\cF|_U)$ to an 
$\cO_X$-algebra $\cJ$, and the $h$ patch to a map 
of $\cO_B$-algebras $\tau\:\cO_X\to\cJ$. 
Since the $\phi_{h,h'}$ do not decrease degrees, 
for each integer $n\geq 0$ the truncated sheaves 
$\prod_{i=0}^n\cS^i(\cF|_U)$ patch to an 
$\cO_X$-algebra quotient $\cJ^n$ of $\cJ$. 
Also, since the $(\phi_{h,h'})_0$ are the identity 
maps, there is a natural map of $\cO_X$-modules 
$\cS^n(\cF)\to\cJ^n$. This map is an isomorphism 
for $n=0$. Also, for each integer $n>0$, the 
$\cO_X$-algebra $\cJ^{n-1}$ is a subquotient of 
$\cJ^n$, and there is a natural exact sequence,
       	$$0\to\cS^n(\cF)\to\cJ^n\to\cJ^{n-1}\to0.$$
We say that $\cJ$ 
is the \emph{sheaf of jets of $D$}, and that 
$\tau\:\cO_X\to\cJ$ is its \emph{Hasse derivation}. 
For each integer $n\geq 0$, we say that $\cJ^n$ 
is the \emph{sheaf of $n$-jets of $D$}, and that the 
induced $\tau_n\:\cO_X\to\cJ^n$ is the 
\emph{n-th order truncated Hasse derivation}.
\end{subsec}

\begin{example}\label{jetsfol}
\setcounter{equation}{0} Let $M$ be a complex 
manifold of complex dimension $m$, 
and $\cL$ a foliation of dimension $n$ 
of $M$. Let $w\:T_M^*\to E$ be the surjection 
associated to $\cL$. Then $w$ induces a 
derivation $D\:\cO_M\to\cF$, where $\cF$ is 
the sheaf of holomorphic sections of $E$. The 
Frobenius conditions imply that $\cF$ is integrable. 
Also, since $w$ is surjective, 
$\cF$ is generated by $D(\cO_M)$. Thus, applying 
Construction~\ref{jets}, 
we have an associated sheaf of $q$-jets $\cJ^q$ on 
$M$ for each integer $q\geq 0$. 
Also we may consider the bundle of $q$-jets 
$J^q_\cL$ of $\cL$. Then $\cJ^q$ is the sheaf 
of holomorphic sections of $J^q_\cL$.

Indeed, for each point $p$ of $M$, there exist 
a neighborhood $X$ of $p$ in 
$M$, and an open embedding of $X$ in 
$\CC^n\times\CC^{m-n}$ whose composition 
$\varphi_X\:X\to\CC^{m-n}$ with the 
second projection defines $\cL$ on $X$. 
The sheaf $\cF|_X$ is the pullback of 
the sheaf of 1-forms on $\CC^n$, and the 
canonical vector fields on $\CC^n$ yield a 
basis $e_1,\dots,e_n$ of $\cF|_X$ such that 
$D=D_1e_1+\cdots+D_ne_n$, where the $D_i$ are 
the pullbacks of these vector fields. Since they 
commute, so do the $D_i$. 

Using the $D$-connection $\gamma$ given in 
Example \ref{loccon}, and the associated extended 
$D$-connection $T$ given in 
Construction \ref{extcon}, we obtain an iterated  
Hasse derivation $h_X$ on $X$ extending $D|_X$, and 
given by Formula~\eqref{hcom} for 
each integer $q>0$ and each local section $a$ of 
$\cF$. Then the truncation in order $q$ of $h_X$ 
has exactly the same form of the canonical Hasse 
derivation of the sheaf of sections of the bundle 
of relative $q$-jets $J^q_{\varphi_X}$. 
So, patching the $h_X$ is compatible with patching 
the $J^q_{\varphi_X}$. The patching of the latter 
yields the bundle of jets $J^q_\cL$, and hence we 
get that 
$\cJ^q$ is the sheaf of sections of $J^q_\cL$.
\end{example}

\end{document}